\documentstyle[12pt]{article}
\begin{document}
\includegraphics{zibheader.eps}
\vspace*{6.5cm}
\begin{center}
{\Large Wolfram Koepf}\vspace*{2cm}\\
\Large{{\bf REDUCE Package for the Indefinite}}\\[2mm]
\Large{{\bf and Definite Summation}}
\end{center}
\vfill
\hrule
\vspace*{3mm}
Technical Report TR 94--9 (November 1994)

\thispagestyle{empty}
\setcounter{page}{0}
\eject
\hoffset -1cm
\voffset -2cm
\footskip=1cm
\parindent=0pt
\title{REDUCE Package for the Indefinite\\
and Definite Summation}
\date{}
\author{Wolfram Koepf\\
	email: {\tt  koepf@zib-berlin.de}
}
\maketitle
\newcommand{\N} {{\rm {\mbox{\protect\makebox[.15em][l]{I}N}}}}

\section*{Abstract}

This article describes the REDUCE package ZEILBERG implemented
by Gregor St\"olting and the author.

The REDUCE package ZEILBERG is a careful implementation of the Gosper%
\footnote{The {\tt sum} package contains also a partial implementation
of the Gosper algorithm.}
and Zeilberger algorithms for indefinite, and definite summation of
hypergeometric terms, respectively. An expression $a_k$ is called a 
{\sl hypergeometric term} (or {\sl closed form}),
if $a_{k}/a_{k-1}$ is a rational function with respect to $k$.
Typical hypergeometric terms are ratios of products of powers, factorials, 
$\Gamma$ function terms, binomial coefficients, and shifted factorials 
(Pochhammer symbols) that are integer-linear in their arguments.

\section{Gosper Algorithm}

The Gosper algorithm \cite{Gos} is a {\sl decision procedure}, that
decides by algebraic calculations whether a given hypergeometric term
$a_k$ has a hypergeometric term antidifference $g_k$, i.\ e.\ 
$g_{k}-g_{k-1}=a_k$, and returns $g_k$ if the procedure is successful, in which
case we call $a_k$ {\sl Gosper-summable}. Otherwise
{\sl no hypergeometric term antidifference exists}. Therefore
if the Gosper algorithm does not return a closed form solution,
it has {\sl proved} that no such solution exists, an information
that may be quite useful and important. 
The Gosper algorithm is the discrete analogue of the Risch algorithm
for integration in terms of elementary functions.

Any antidifference is uniquely determined up to a constant, and is
denoted by
\[
g_k=\sum\nolimits_k a_k
\;.
\]
Finding $g_k$ given $a_k$ is called {\sl indefinite summation}.
The antidifference operator $\Sigma$ is the inverse of the downward 
difference operator $\nabla a_k=a_{k}-a_{k-1}$. There is an analogous
summation theory corresponding to the upward difference operator 
$\Delta a_k=a_{k+1}-a_k$.

In case, an antidifference $g_k$ of $a_k$ is known, any sum
\[
\sum_{k=m}^{n} a_k=g_{n}-g_{m-1}
\]
can be easily calculated by an evaluation of $g$ at the boundary points
like in the integration case. Note, however, that the sum
\begin{equation}
\sum_{k=0}^n {{n}\choose{k}}
\label{eq:nchoosek}
\end{equation}
e.\ g.\
is not of this type as the summand ${{n}\choose{k}}$ depends on the upper
boundary point $n$ explicitly. This is an example of a definite sum
that we consider in the next section.

Our package supports the input of powers ({\tt a\verb+^+k)},
factorials ({\tt factorial(k)}),
$\Gamma$ function terms ({\tt gamma(a)}), binomial coefficients
({\tt binomial(n,k)}), shifted factorials 
({\tt pochhammer(a,k)$=a(a+1)\cdots(a+k-1)=\Gamma (a+k)/\Gamma (a)$}), and
partially products ({\tt prod(f,k,k1,k2)}).
It takes care of the necessary simplifications, and therefore 
(in principle) provides you with the solution of the decision problem
as long as the memory or time requirements are not too high for the
computer used.

\section{Zeilberger Algorithm}

The (fast) Zeilberger algorithm \cite{Zei2}--\cite{Zei3}
deals with the {\sl definite summation} of 
hypergeometric terms. Zeilberger's paradigm is to find (and return)
a linear homogeneous recurrence equation with polynomial coefficients
(called {\sl holonomic equation}) for an {\sl infinite sum}
\[
s(n)=\sum_{k=-\infty}^{\infty} f(n,k)
\;,
\]
the summation to be understood over all integers $k$,
if $f(n,k)$ is a hypergeometric term with respect to both $k$ and $n$.
The existence of a holonomic recurrence equation for $s(n)$ is then
generally guaranteed.

If one is lucky, and the resulting recurrence equation is of first order
\[
p(n)\,s(n-1)+q(n)\,s(n)=0
\quad\quad(p,q\;\mbox{polynomials})
\;,
\]
$s(n)$ turns out to be a hypergeometric term, and a closed form solution
can be easily established using a suitable initial value, and is
represented by a ratio of Pochhammer or $\Gamma$ function terms if the
polynomials $p$, and $q$ can be factored.

Zeilberger's algorithm does not guarantee to find the holonomic equation
of lowest order, but often it does. 

If the resulting recurrence equation has order larger than one, 
this information can be used for identification purposes:
Any other expression satisfying the same recurrence equation, and the same
initial values, represents the same function.

Note that a {\sl definite sum} $\sum\limits_{k=m_1}^{m_2} f(n,k)$ is an
infinite sum if $f(n,k)=0$ for $k<m_1$ and $k>m_2$.
This is often the case, an example of which is the sum (\ref{eq:nchoosek})
considered above, for which the hypergeometric recurrence equation
$2 s(n-1) - s(n) = 0$ is generated by Zeilberger's algorithm, leading 
to the closed form solution $s(n)=2^n$.

Definite summation is trivial if the corresponding indefinite sum
is Gosper-summable analogously to the fact that definite integration
is trivial as soon as an elementary antiderivative is known.  If this is 
not the case, the situation is much more difficult, and it is therefore
quite remarkable and non-obvious
that Zeilberger's method is just a clever application of Gosper's algorithm.

Our implementation is mainly based on \cite{Koornwinder}. Many more
examples can be found in \cite{PS}, \cite{Strehl2}, \cite{Wil1},
and \cite{Wilf} most of which are contained in the test file 
{\tt zeilberg.tst}.

\section{REDUCE operator {\tt GOSPER}}
The ZEILBERG package must be loaded by:

{\small
\begin{verbatim}
1: load zeilberg;
\end{verbatim}
}\noindent
The {\tt gosper} operator is an implementation of the Gosper algorithm.
\begin{itemize}
\item
{\tt gosper(f,k)} determines a closed
form antidifference. If it does not return a closed form solution, then
a closed form solution does not exist.
\item
{\tt gosper(f,k,m,n)} determines 
\[
\sum_{k=m}^n a_k
\]
using Gosper's algorithm. This is only successful if Gosper's algorithm applies.
\end{itemize}

Example:

{\small
\begin{verbatim}
2: gosper((-1)^(k+1)*(4*k+1)*factorial(2*k)/ 
   (factorial(k)*4^k*(2*k-1)*factorial(k+1)),k);

             k
     - ( - 1) *factorial(2*k)
----------------------------------
  k
 4 *factorial(k + 1)*factorial(k)
\end{verbatim}
}\noindent
This solves a problem given in SIAM Review (\cite{SR}, Problem 94--2) 
where it was asked to determine the infinite sum
\[
S=\lim_{n\rightarrow\infty} S_n
\;,
\quad\quad\quad
S_n=\sum_{k=1}^n
\frac{(-1)^{k+1}(4k+1)(2k-1)!!}{2^k(2k-1)(k+1)!}
\;,
\]
($(2k-1)!!=1\cdot 3 \cdots (2k-1)=\frac{(2k)!}{2^k\,k!}$).
The above calculation shows that the summand is Gosper-summable,
and the limit $S=1$ is easily established using Stirling's formula.

The implementation solves further deep and difficult problems some examples of
which are:%
\footnote{Note that REDUCE Version 3.5 gives the output in
terms of $\Gamma$ functions.}

{\small
\begin{verbatim}
3:  gosper(sub(n=n+1,binomial(n,k)^2/binomial(2*n,n))-
    binomial(n,k)^2/binomial(2*n,n),k);

                   2                                                          2
((binomial(n + 1,k) *binomial(2*n,n) - binomial(2*(n + 1),n + 1)*binomial(n,k) )

                             2
 *(2*k - 3*n - 1)*(k - n - 1) )/(

                                                     2
   (2*(2*(n + 1) - k)*(2*n + 1)*k - (3*n + 1)*(n + 1) )

   *binomial(2*(n + 1),n + 1)*binomial(2*n,n))

4: gosper(binomial(k,n),k);

 (k + 1)*binomial(k,n)
-----------------------
         n + 1
5: gosper((-25+15*k+18*k^2-2*k^3-k^4)/
   (-23+479*k+613*k^2+137*k^3+53*k^4+5*k^5+k^6),k);

          2
    - (2*k  - 15*k + 8)*k
----------------------------
      3      2
 23*(k  + 4*k  + 27*k + 23)
\end{verbatim}
}\noindent
The Gosper algorithm is not able to give antidifferences depending
on the harmonic numbers
\[
H_k:=\sum_{j=1}^k\frac{1}{j}
\;,
\]
e.\ g.\ $\sum_k H_k=(k+1)(H_{k+1}-1)$, but, is able to give a proof, instead,
for the fact that $H_k$ does not possess a closed form evaluation:

{\small
\begin{verbatim}
6: gosper(1/k,k);

***** Gosper algorithm: no closed form solution exists
\end{verbatim}
}\noindent
The following code gives the solution to a summation problem proposed in 
Gosper's original paper \cite{Gos}. Let
\[
f_k=\prod_{j=1}^k (a+b\,j+c\,j^2)
\quad\quad\mbox{and}\quad\quad
g_k=\prod_{j=1}^k (e+b\,j+c\,j^2)
\;.
\]
Then a closed form solution for
\[
\sum\nolimits_k\frac{f_{k-1}}{g_{k}}
\]
is found by the definitions

{\small
\begin{verbatim}
7: operator ff,gg$

8: let {ff(~k+~m) => ff(k+m-1)*(c*(k+m)^2+b*(k+m)+a) when (fixp(m) and m>0),
   ff(~k+~m) => ff(k+m+1)/(c*(k+m+1)^2+b*(k+m+1)+a) when (fixp(m) and m<0)}$

9: let {gg(~k+~m) => gg(k+m-1)*(c*(k+m)^2+b*(k+m)+e) when (fixp(m) and m>0),
   gg(~k+~m) => gg(k+m+1)/(c*(k+m+1)^2+b*(k+m+1)+e) when (fixp(m) and m<0)}$
\end{verbatim}
}\noindent
and the calculation

{\small
\begin{verbatim}
10: gosper(ff(k-1)/gg(k),k);

     ff(k)
---------------
 (a - e)*gg(k)

11: clear ff,gg$
\end{verbatim}
}\noindent
Similarly closed form solutions of $\sum\nolimits_k\frac{f_{k-m}}{g_{k}}$
for positive integers $m$ can be obtained, as well as of
$\sum_k\frac{f_{k-1}}{g_{k}}$ for
\[
f_k=\prod_{j=1}^k (a+b\,j+c\,j^2+d\,j^3)
\quad\quad\mbox{and}\quad\quad
g_k=\prod_{j=1}^k (e+b\,j+c\,j^2+d\,j^3)
\]
and for analogous expressions of higher degree polynomials.

\section{REDUCE operator {\tt SUMRECURSION}}
The {\tt sumrecursion} operator is an implementation of the (fast)
Zeilberger algorithm.
\begin{itemize}
\item
{\tt sumrecursion(f,k,n)} determines a holonomic recurrence equation
for $\sum\limits_{k=-\infty}^\infty f(n,k)$ with respect to $n$.
\item
{\tt sumrecursion(f,k,n,j)} $(j\in\N)$ 
searches only for a holonomic recurrence equation of order $j$.
\end{itemize}
A simple example deals with Equation (\ref{eq:nchoosek})%
\footnote{Note that with REDUCE Version 3.5 we use the global operator 
{\tt summ} instead of {\tt sum} to denote the sum.} 

{\small
\begin{verbatim}
12: sumrecursion(binomial(n,k),k,n);

2*sum(n - 1) - sum(n)
\end{verbatim}
}\noindent
The whole {\sl hypergeometric database} of the {\sl
Vandermonde, Gau{\ss}, Kummer, Saalsch\"utz, Dixon, Clausen} and {\sl Dougall
identities} (see \cite{Wilf}), and many more identities,
can be obtained using {\tt sumrecursion}. 
As examples, we consider the difficult cases of Clausen and Dougall:%
\footnote{Note that the latter may need a large amount of computing time.}

{\small
\begin{verbatim}
13: summand:=factorial(a+k-1)*factorial(b+k-1)/
    (factorial(k)*factorial(-1/2+a+b+k))*
    factorial(a+n-k-1)*factorial(b+n-k-1)/
    (factorial(n-k)*factorial(-1/2+a+b+n-k))$

14: sumrecursion(summand,k,n);

(2*a + 2*b + 2*n - 1)*(2*a + 2*b + n - 1)*sum(n)*n

 - 2*(2*a + n - 1)*(a + b + n - 1)*(2*b + n - 1)*sum(n - 1)

15: summand:=pochhammer(d,k)*pochhammer(1+d/2,k)*
    pochhammer(d+b-a,k)*pochhammer(d+c-a,k)*
    pochhammer(1+a-b-c,k)*pochhammer(n+a,k)*
    pochhammer(-n,k)/(factorial(k)*pochhammer(d/2,k)*
    pochhammer(1+a-b,k)*pochhammer(1+a-c,k)*
    pochhammer(b+c+d-a,k)*pochhammer(1+d-a-n,k)*pochhammer(1+d+n,k))$

16: sumrecursion(summand,k,n);

(2*a - b - c - d + n)*(b + n - 1)*(c + n - 1)*(d + n)*sum(n - 1)

 + (a - b - c - d - n + 1)*(a - b + n)*(a - c + n)*(a - d + n - 1)*sum(n)
\end{verbatim}
}\noindent
corresponding to the statements
\[
_4 F_3\left.
\!\!
\left(
\!\!\!\!
\begin{array}{c}
\multicolumn{1}{c}{\begin{array}{c}
a\;, b\;, 1/2-a-b-n\;, -n
\end{array}}\\[1mm]
\multicolumn{1}{c}{\begin{array}{c}
1/2+a+b \;, 1-a-n\;, 1-b-n
            \end{array}}\end{array}
\!\!\!\!
\right| 1\right)
=\frac{(2a)_n\,(a+b)_n\,(2b)_n}
{(2a+2b)_n\,(a)_n\,(b)_n}
\]
and
\[
_7 F_6\left.
\!\!
\left(
\!\!\!\!
\begin{array}{c}
\multicolumn{1}{c}{\begin{array}{c} 
d\;, 1+d/2\;, d+b-a\;, d+c-a\;, 1+a-b-c\;, n+a\;, -n 
\end{array}}\\[1mm]
\multicolumn{1}{c}{\begin{array}{c} 
 d/2\;, 1+a-b\;, 1+a-c\;, b+c+d-a \;, 1+d-a-n\;, 1+d+n
            \end{array}}\end{array}
\!\!\!\!
\right| 1\right)
\]
\[
=\frac{(d+1)_n\,(b)_n\,(c)_n\,(1+2\,a-b-c-d)_n}
{(a-d)_n\,(1+a-b)_n\,(1+a-c)_n\,(b+c+d-a)_n}
\]
(compare next section), respectively.

Other applications of the Zeilberger algorithm are connected with
the verification of identities. To prove the identity
\[
\sum_{k=0}^n
{{n}\choose{k}}^3
=
\sum_{k=0}^n
{{n}\choose{k}}^2 {{2k}\choose{n}}
\;,
\]
e.\ g., we may prove that both sums satisfy the same recurrence equation

{\small
\begin{verbatim}
17: sumrecursion(binomial(n,k)^3,k,n);

         2                        2       2
8*(n - 1) *sum(n - 2) - sum(n)*n  + (7*n  - 7*n + 2)*sum(n - 1)

18: sumrecursion(binomial(n,k)^2*binomial(2*k,n),k,n);

         2                        2       2
8*(n - 1) *sum(n - 2) - sum(n)*n  + (7*n  - 7*n + 2)*sum(n - 1)
\end{verbatim}
}\noindent
and finally check the initial conditions:

{\small
\begin{verbatim}
19: sum(sub(n=0,binomial(n,k)^3),k,0,0);

1

20: sum(sub(n=0,binomial(n,k)^2*binomial(2*k,n)),k,0,0);

1

21: sum(sub(n=1,binomial(n,k)^3),k,0,1);

2

22: sum(sub(n=1,binomial(n,k)^2*binomial(2*k,n)),k,0,1);

2
\end{verbatim}
}\noindent

\section{REDUCE operator {\tt HYPERRECURSION}}
Sums for which the Zeilberger algorithm applies, in general are
special cases of the {\sl generalized hypergeometric function}
\[
_{p}F_{q}\left.\left(\begin{array}{cccc}
a_{1},&a_{2},&\cdots,&a_{p}\\
b_{1},&b_{2},&\cdots,&b_{q}\\
            \end{array}\right| x\right)
:=
\sum_{k=0}^\infty \frac
{(a_{1})_{k}\cdot(a_{2})_{k}\cdots(a_{p})_{k}}
{(b_{1})_{k}\cdot(b_{2})_{k}\cdots(b_{q})_{k}\,k!}x^{k}
\label{eq:coefficientformula}
\]
with upper parameters $\{a_{1}, a_{2}, \ldots, a_{p}\}$, and lower
parameters $\{b_{1}, b_{2}, \ldots, b_{q}\}$. If a recursion for a
generalized hypergeometric function is to be established, you can use
the following REDUCE operator:
\begin{itemize}
\item
{\tt hyperrecursion(upper,lower,x,n)} determines a holonomic recurrence 
equation with respect to $n$ for 
$_{p}F_{q}\left.\left(\begin{array}{cccc}
a_{1},&a_{2},&\cdots,&a_{p}\\
b_{1},&b_{2},&\cdots,&b_{q}\\
            \end{array}\right| x\right)
$, where {\tt upper}$=\{a_{1}, a_{2}, \ldots, a_{p}\}$
is the list of upper parameters, and 
{\tt lower}$=\{b_{1}, b_{2}, \ldots, b_{q}\}$
is the list of lower parameters depending on $n$.
\item
{\tt hyperrecursion(upper,lower,x,n,j)} $(j\in\N)$
searches only for a holonomic recurrence equation of order $j$.
\end{itemize}
Therefore

{\small
\begin{verbatim}
23: hyperrecursion({-n,b},{c},1,n);

(n - 1 + c - b)*sum(n - 1) - (n - 1 + c)*sum(n)
\end{verbatim}
}\noindent
establishes the Vandermonde identity
\[
_2 F_1\left.
\!\!
\left(
\!\!\!\!
\begin{array}{c}
\multicolumn{1}{c}{\begin{array}{cc} -n\;, & b \end{array}}\\[1mm]
\multicolumn{1}{c}{ c}
            \end{array}
\!\!\!\!
\right| 1\right)
=\frac{(c-b)_n}{(c)_n}
\;,
\]
whereas

{\small
\begin{verbatim}
24: hyperrecursion({d,1+d/2,d+b-a,d+c-a,1+a-b-c,n+a,-n},
                   {d/2,1+a-b,1+a-c,b+c+d-a,1+d-a-n,1+d+n},1,n);

(2*a - b - c - d + n)*(b + n - 1)*(c + n - 1)*(d + n)*sum(n - 1)

 + (a - b - c - d - n + 1)*(a - b + n)*(a - c + n)*(a - d + n - 1)*sum(n)
\end{verbatim}
}\noindent
proves Dougall's identity, again.

If a hypergeometric expression is given in hypergeometric notation, then
the use of\linebreak {\tt hyperrecursion} is more natural than the use of
{\tt sumrecursion}.

Moreover you may use the REDUCE operator
\begin{itemize}
\item
{\tt hyperterm(upper,lower,x,k)} that yields the hypergeometric term
\[
\frac
{(a_{1})_{k}\cdot(a_{2})_{k}\cdots(a_{p})_{k}}
{(b_{1})_{k}\cdot(b_{2})_{k}\cdots(b_{q})_{k}\,k!}x^{k}
\]
with upper parameters {\tt upper}$=\{a_{1}, a_{2}, \ldots, a_{p}\}$,
and lower parameters {\tt lower}$=\linebreak \{b_{1}, b_{2}, \ldots, b_{q}\}$
\end{itemize}
in connection with hypergeometric terms.

The operator {\tt sumrecursion} can also be used to 
obtain three-term recurrence equations for systems of orthogonal polynomials
with the aid of known hypergeometric representations. By
(\cite{NSU}, (2.7.11a)), the discrete Krawtchouk polynomials $k_n^{(p)}(x,N)$
have the hypergeometric representation
\[
k_n^{(p)}(x,N)=
(-1)^n\,p^n\,{{N}\choose{n}}\;
_2 F_1\left.
\!\!
\left(
\!\!\!\!
\begin{array}{c}
\multicolumn{1}{c}{\begin{array}{cc} -n\;, & -x \end{array}}\\[1mm]
\multicolumn{1}{c}{ -N}
            \end{array}
\!\!\!\!
\right| \frac{1}{p}\right)
\;,
\]
and therefore we declare

{\small
\begin{verbatim}
25: krawtchoukterm:=(-1)^n*p^n*binomial(NN,n)*hyperterm({-n,-x},{-NN},1/p,k)$
\end{verbatim}
}\noindent
and get the three three-term recurrence equations

{\small
\begin{verbatim}
26: sumrecursion(krawtchoukterm,k,n);

(x + 1 - 2*p - nn*p + (2*p - 1)*n)*sum(n - 1)

 - ((n - nn - 2)*(p - 1)*sum(n - 2)*p + sum(n)*n)

27: sumrecursion(krawtchoukterm,k,x);

 - ((x - 1 + nn*p - n - 2*(x - 1)*p)*sum(x - 1) + (x - 1 - nn)*sum(x)*p

     + (p - 1)*(x - 1)*sum(x - 2))

28: sumrecursion(krawtchoukterm,k,NN);

(x + 1 + n + (p - 2)*nn)*sum(nn - 1)

 - ((x + 1 - nn)*sum(nn - 2) - (n - nn)*(p - 1)*sum(nn))
\end{verbatim}
}\noindent
with respect to the parameters $n$, $x$, and $N$ respectively.

\section{Simplification Operators}
For the decision that an expression $a_k$ is a hypergeometric term, it is 
necessary to find out whether or not $a_{k}/a_{k-1}$ is a rational
function with respect to $k$. For the purpose to decide
whether or not an expression involving powers, factorials,
$\Gamma$ function terms, 
binomial coefficients, and Pochhammer symbols is a hypergeometric term,
the following simplification operators can be used:
\begin{itemize}
\item
{\tt simplify\verb+_+gamma(f)} simplifies an expression {\tt f} involving 
only rational, powers and\linebreak
$\Gamma$ function terms according to a recursive 
application of the simplification rule $\Gamma\:(a+1)=a\,\Gamma\:(a)$ 
to the expression tree. Since all $\Gamma$ arguments with integer difference
are transformed, this gives a decision procedure for rationality
for integer-linear $\Gamma$ term product ratios.
\item
{\tt simplify\verb+_+combinatorial(f)} simplifies an expression {\tt f} 
involving powers, factorials,\linebreak
$\Gamma$ function terms,
binomial coefficients, and Pochhammer symbols by converting
factorials, binomial coefficients, and Poch\-hammer symbols into 
$\Gamma$ function terms, and applying {\tt simplify\verb+_+gamma} to its
result. If the output is not rational,
it is given in terms of $\Gamma$ functions. If you prefer factorials
you may use
\item
{\tt gammatofactorial} (rule) converting $\Gamma$ function terms in
factorials using $\Gamma\:(x)\rightarrow (x-1)!$.
\end{itemize}
The use of {\tt simplify\verb+_+combinatorial(f)} is a safe way to
decide the rationality for any ratio of products of powers, factorials,
$\Gamma$ function terms, binomial coefficients, and Pochhammer symbols.

\pagebreak
Example:

{\small
\begin{verbatim}
29: simplify_combinatorial(sub(k=k+1,krawtchoukterm)/krawtchoukterm);

  (k - n)*(k - x)
--------------------
 (k - nn)*(k + 1)*p
\end{verbatim}
}\noindent
From this calculation, we see again that the upper parameters of
the hypergeometric representation of the Krawtchouk polynomials are given by
$\{-n,-x\}$, its lower parameter is $\{-n\}$, and the argument of the
hypergeometric function is $1/p$.

Another example is

{\small
\begin{verbatim}
30: simplify_combinatorial(binomial(n,k));

         gamma(n + 1)
-------------------------------
 gamma(k + 1)*gamma(n + 1 - k)

31: ws where gammatofactorial;

            factorial(n)
------------------------------------
 factorial( - (k - n))*factorial(k)
\end{verbatim}
}\noindent

\section{Tracing}
If you set 

{\small
\begin{verbatim}
32: on zb_trace;
\end{verbatim}
}\noindent
tracing is enabled, and you get intermediate results, see \cite{Koornwinder}. 

Example for the Gosper algorithm:

{\small
\begin{verbatim}
33: gosper(pochhammer(k-n,n),k);

                 k - 1
a(k)/a(k-1):= -----------
               k - n - 1

Gosper algorithm applicable

p:= 1

q:= k - 1

r:=  - (n + 1 - k)

degreebound := 0

       1
f:= -------
     n + 1

Gosper algorithm successful

 pochhammer(k - n,n)*k
-----------------------
         n + 1
\end{verbatim}
}\noindent
Example for the Zeilberger algorithm:

{\small
\begin{verbatim}
34: sumrecursion(binomial(n,k)^2,k,n);

                          2
                         n
F(n,k)/F(n-1,k):= -----------------
                    2            2
                   k  - 2*k*n + n

                    2                  2
                   k  - 2*k*n - 2*k + n  + 2*n + 1
F(n,k)/F(n,k-1):= ---------------------------------
                                  2
                                 k

Zeilberger algorithm applicable

applying Zeilberger algorithm for order:= 1

                 2                                    2    2
p:= zb_sigma(1)*k  - 2*zb_sigma(1)*k*n + zb_sigma(1)*n  + n

     2                  2
q:= k  - 2*k*n - 2*k + n  + 2*n + 1

     2
r:= k

degreebound := 1

     2*k - 3*n + 2
f:= ---------------
           n

           2        2        2              3      2
      - 4*k *n + 2*k  + 8*k*n  - 4*k*n - 3*n  + 2*n
p:= -------------------------------------------------
                            n

Zeilberger algorithm successful

4*sum(n - 1)*n - 2*sum(n - 1) - sum(n)*n

35: off zb_trace;
\end{verbatim}
}\noindent

\section{Global Variables and Switches}
The following global variables and switches can be used in connection with
the {\tt ZEILBERG} package:
\begin{itemize}
\item
{\tt zb\verb+_+trace}, switch; default setting {\tt on}. 
Turns tracing
on and off.
\item
{\tt zb\verb+_+direction}, variable; settings: {\tt down}, {\tt up};
default setting {\tt down}. 

In the case of the Gosper algorithm, either a downward or a forward
antidifference is calculated, i.\ e., {\tt gosper} finds $g_k$ with either
\[
a_k=g_k-g_{k-1}
\quad\quad\mbox{or}\quad\quad
a_k=g_{k+1}-g_{k},
\]
respectively.

In the case of the Zeilberger algorithm, either a downward or an upward
recurrence equation is given. Example:

{\small
\begin{verbatim}
36: zb_direction:=up$

37: sumrecursion(binomial(n,k)^2,k,n);

sum(n + 1)*n + sum(n + 1) - 4*sum(n)*n - 2*sum(n)

38: zb_direction:=down$
\end{verbatim}
}\noindent
\item
{\tt zb\verb+_+order}, variable; settings: any nonnegative integer;
default setting~{\tt 5}. 
Gives the maximal order for the recurrence
equation that {\tt sumrecursion} searches for.
\item
{\tt zb\verb+_+factor}, switch; default setting {\tt on}.
If {\tt off}, the factorization of the output usually producing nicer results
is suppressed.
\item
{\tt zb\verb+_+proof}, switch; default setting {\tt off}. If {\tt on},
then several intermediate results are stored in global variables:
\item
{\tt gosper\verb+_+representation}, variable; default setting {\tt nil}.

If a {\tt gosper} command is issued, and if the Gosper algorithm is applicable,
then the variable {\tt gosper\verb+_+representation} is set to the
list of polynomials (with respect to $k$) {\tt \{p,q,r,f\}} 
corresponding to the representation
\[
\frac{a_k}{a_{k-1}}=\frac{p_k}{p_{k-1}}\,\frac{q_k}{r_k}
\;,
\quad\quad\quad
g_k=\frac{q_{k+1}}{p_k}\,f_k\,a_k
\;,
\]
see \cite{Gos}. Examples:

{\small
\begin{verbatim}
39: on zb_proof;

40: gosper(k*factorial(k),k);

(k + 1)*factorial(k)

41: gosper_representation;

{k,k,1,1}

42: gosper(1/(k+1)*binomial(2*k,k)/(n-k+1)*binomial(2*n-2*k,n-k),k);

 (2*k - n + 1)*(2*k + 1)*binomial( - 2*(k - n), - (k - n))*binomial(2*k,k)
---------------------------------------------------------------------------
                          (k + 1)*(n + 2)*(n + 1)

43: gosper_representation;

{1,

 (2*k - 1)*(k - n - 2),

 (2*k - 2*n - 1)*(k + 1),

    n - 1 - 2*k
 -----------------}
  (n + 2)*(n + 1)
\end{verbatim}
}\noindent
\item
{\tt zeilberger\verb+_+representation}, variable; default setting {\tt nil}.

If a {\tt sumrecursion} command is issued, and if the Zeilberger
algorithm is successful, then the variable 
{\tt zeilberger\verb+_+representation} is set to the final Gosper
representation used, see \cite{Koornwinder}.
\end{itemize}

\pagebreak
\section{Messages}

The following messages may occur:
\begin{itemize}
\item
{\tt ***** Gosper algorithm:\ no closed form solution exists}

Example input:

{\tt gosper(factorial(k),k)}.
\item
{\tt ***** Gosper algorithm not applicable}

Example input:

{\tt gosper(factorial(k/2),k)}.

The term ratio $a_k/a_{k-1}$ is not rational.
\item
{\tt ***** illegal number of arguments}

Example input:

{\tt gosper(k)}.
\item
{\tt ***** Zeilberger algorithm fails.\ Enlarge zb\verb+_+order}

Example input:

{\tt sumrecursion(binomial(n,k)*binomial(6*k,n),k,n)}

For  this example a setting {\tt zb\verb+_+order:=6} is needed.
\item
{\tt ***** Zeilberger algorithm not applicable}

Example input:

{\tt sumrecursion(binomial(n/2,k),k,n)}

One of the term ratios $f(n,k)/f(n-1,k)$ or $f(n,k)/f(n,k-1)$
is not rational.
\end{itemize}

\section*{Acknowledgement}
I like to thank Gregor St\"olting for his careful implementation,
Winfried Neun for his assistance concerning internal details of REDUCE,
and Prof.\ Peter Deuflhard for his support of my research.

\pagebreak

\end{document}